# SOBOL' INDICES FOR PROBLEMS DEFINED IN NON-RECTANGULAR DOMAINS

*S. Kucherenko, O.V. Klymenko, N. Shah*
*Imperial College London, London, SW7 2AZ, UK*
e-mail: s.kucherenko@imperial.ac.uk

**Abstract.** A novel theoretical and numerical framework for the estimation of Sobol' sensitivity indices for models in which inputs are confined to a non-rectangular domain (e.g., in presence of inequality constraints) is developed. Two numerical methods, namely the quadrature integration method which may be very efficient for problems of low and medium dimensionality and the MC/QMC estimators based on the acceptance-rejection sampling method are proposed for the numerical estimation of Sobol' sensitivity indices. Several model test functions with constraints are considered for which analytical solutions for Sobol' sensitivity indices were found. These solutions were used as benchmarks for verifying numerical estimates. The method is shown to be general and efficient.

## 1. Introduction

High complexity of models in physics, chemistry, environmental studies, reliability and structural analysis and other fields results in the increased uncertainty in model parameters and model structures. Uncertainty in inputs is reflected in uncertainty of model outputs and predictions. Uncertainty and sensitivity analysis has been recognized as an essential part of model applications. Global sensitivity analysis (GSA) offers a comprehensive approach to model analysis by quantifying how the uncertainty in model output is apportioned to the uncertainty in model inputs [1,2]. Unlike local sensitivity analysis, GSA estimates the effect of varying a given input (or set of inputs) while all other inputs are varied as well, thus providing a measure of interactions among variables. GSA is used to identify key parameters whose uncertainty most affects the output. This information then can be used to rank variables, fix unessential variables and thus decrease problem dimensionality. The variance-based method based on Sobol' sensitivity indices has become very popular among practitioners due to its efficiency and ease of interpretation [3,4]. Most of the developed techniques for GSA are designed under the hypothesis that inputs are independent. However, in many cases there are dependencies among inputs, which may have significant impact on the importance results.

There have been a number of attempts to extend GSA to models with dependent inputs. We present a brief survey of only the most recent and easy-to-follow developments. Xu and Gertner [5] suggested to split the contribution of an individual input to the uncertainty of the model output into two components: the correlated contribution and the uncorrelated one. A regression-based method was used for estimating the correlated and uncorrelated contributions of the inputs. The approach developed originally only for linear models was later extended to nonlinear models [6]. Li *et al* [7] proposed a generalization of the ANOVA-HDMR decomposition by including covariances to the model variance. They distinguished between structural and correlative contributions of a given input. This method presents some critical issues such as non-uniqueness of the functional decomposition and hence difficulties in interpreting the results. Da Veiga *et al*. [8] used the estimator of the main effect index firstly proposed in [9] and further developed in [10]. They suggested a new procedure based on a heteroscedastic regression model and local polynomial metamodel.

A novel approach for the estimation of Sobol' sensitivity indices for models with dependent variables using generalizations of direct Sobol' formulas was developed in Kucherenko *et al.* [11]. Both the main effect and total sensitivity indices were derived as generalizations of Sobol' sensitivity indices. Formulas and Monte Carlo (MC) numerical estimators similar to Sobol' formulas were proposed. A Gaussian copula-based approach was used for sampling from arbitrary multivariate probability distributions. The generalization does not involve the use of surrogate models, data-fitting procedures or orthogonalization of the input factor space.

Mara and Tarantola [12] introduced a set of sensitivity indices which relies on orthogonalisation of correlated inputs. The computation of sensitivity indices was performed using a parametric method, specifically the polynomial chaos expansion. Mara *et al*. [13] extended the development of sensitivity indices suggested in [12] by proposing two sampling strategies for non-parametric, numerical estimation using the Rosenblatt transformation (RT) [14]. RT, and hence values of sensitivity indices, is not unique: for a model with *n* inputs there are *n!* possibilities corresponding to all possible permutations of the elements of the input vector $x = (x_1,...,x_n)$. The authors considered only the RT obtained after circularly reordering the set $(x_1,...,x_n)$, resulting in *n* RT transformations. The authors also established the link with the indices proposed by Kucherenko et al. [11]. W.Hao *et al.* [15] suggested a detailed interpretation of indices proposed in [11]. Taking a quadratic polynomial model with normal inputs as an illustration, they derived explicit expressions for sensitivity indices and considered contributions to the values of sensitivity indices from all components.

In this work we propose GSA formulations for an even wider class of problems with dependent variables, namely for models involving inequality constraints (which naturally leads to the term 'constrained GSA' or cGSA). Such constraints impose structural dependences between model variables in addition to potential correlations between them. This implies that the parameter space may no longer be considered as an *n*-dimensional hypercube or hyperrectangle as considered within GSA so far, but may assume any shapes (including disconnected domains) depending on the number and nature of constraints. This class of problems encompasses a wide range of situations encountered in the natural sciences, engineering, design, economics and finances where model variables are subject to certain limitations imposed e.g. by conservation laws, geometry, costs, quality constraints etc. An important particular case within this paradigm corresponds to imposing a minimum (maximum) threshold for the model output, i.e., $f(x_1,...,x_n) \geq f_{\min}$, in which case the constraint function can be defined as $g(x_1,...,x_n) = f(x_1,...,x_n) - f_{\min}$ and the corresponding constraint as $g(x_1,...,x_n) \geq 0$.

The development of efficient computational approaches for cGSA is challenging because of potentially arbitrary shape of the feasible domain of the model variables' variation, thus requiring the development of special MC or quasi MC (QMC) sampling techniques and approaches for computing sensitivity indices. This becomes especially difficult for models of high dimensionality. Another challenge consists in analysing and interpreting model variance and sensitivity indices in such domains, since in this case sensitivity indices will carry structural information imposed by model constraints. The interpretation of sensitivity indices in such circumstances is crucial for the efficient design of experiments as well as for potential model reduction by fixing or eliminating nonessential variables.

In this paper we have developed and tested several approaches for the numerical estimation of main effect and total sensitivity indices in the cGSA setting. It is organized as follows: The next Section presents formulas for the main effect and total sensitivity indices for models with dependent variables and acceptance–rejection method which can be used to avoid sampling from conditional distributions. Section 3 considers numerical aspects of the method, namely it presents general MC estimators, acceptance-rejection estimation of Sobol' sensitivity indices using grid quadrature formulas and MC estimators. Section 4 presents the application of the proposed methodology to three test cases. Numerical and analytical results are compared and the convergence rates of the MC and QMC methods are discussed. Finally, conclusions are presented in the last Section.

## 2. Variance-based sensitivity indices for models with dependent inputs

### 2.1. Main and total effect formula for inputs defined in non-rectangular domain $\Omega^n$

Consider a model function $f(x_1,...,x_n)$ defined in $R^n$. The input vector $x = (x_1,...,x_n)$ consists of real-valued random variables with a continuous joint probability distribution function (pdf) $p(x_1,...,x_n)$. It is assumed that $f(x_1,...,x_n)$ has a finite variance. Consider an arbitrary subset of the variables $y = (x_{i_1},...,x_{i_s})$, $1 \leq s < n$ and a complementary subset $z = (x_{i_{s+1}},...,x_{i_n})$, so that $x = (y, z)$.

The ratio

$$S_y = \frac{D_y[E_z(f(y,\bar{z}))]}{D} \quad (1)$$

is known as the main effect index of the subset $y$. The quantity

$$S_y^T = \frac{E_z[D_y(f(\bar{y},z))]}{D} \quad (2)$$

is known as the total effect index of the subset $y$. Collectively $S_y$, $S_y^T$ are known as Sobol' indices. In this paper our notations follow closely those of Kucherenko et al. [11].

$E_z(f(y,\bar{z}))$ in Eq. (1) denotes a conditional expectation with respect to $z$ with $y$ being fixed and $D_y(f(\bar{y},z))$ in (2) is a conditional variance with respect to $y$ with $z$ being fixed. We use notations $z$ and $\bar{z}$ (correspondingly $y$ and $\bar{y}$) to distinguish a random vector $z$ (correspondingly $y$) generated from a joint pdf $p(y,z)$ and a random vector $\bar{z}$ (correspondingly $\bar{y}$) generated from a conditional pdf $p(y,\bar{z}|y)$ (correspondingly $p(\bar{y},z|z)$).

A formula for the main effect index can be explicitly written as

$$S_y = \frac{1}{D}\left[\int_{R^s} p(y)dy\left[\int_{R^{n-s}} f(y,\bar{z})p(y,\bar{z}|y)d\bar{z}\right]^2 - f_0^2\right]. \quad (3)$$

Here $f_0 = E(f(x))$, $p(y)$ is the marginal pdf of subset of inputs $y$. This is the so-called double loop formula (DL) [1]. It can be transformed into a different formula which in some cases can be more efficient [11]:

$$S_y = \frac{1}{D}\left[\int_{R^n} f(y',z')p(y',z')dy'dz'\left[\int_{R^{n-s}} f(y',\hat{z})p(y',\hat{z}|y')d\hat{z} - \int_{R^n} f(y,z)p(y,z)dydz\right]\right]. \quad (4)$$

The following notation is used: $z$, $z'$ are two different random vectors generated from the joint pdf $p(y,z)$, a random vector $\hat{z}$ is generated from a conditional pdf $p(y',\hat{z}|y')$.

A formula for computing the total effect index $S_y^T$ (2) can be written as

$$S_y^T = \frac{1}{2D}\int_{R^{n+s}} [f(y,z) - f(\bar{y},z)]^2 p(y,z)p(\bar{y},z|z)dyd\bar{y}dz. \quad (5)$$

Consider a more general situation when the variation of model inputs is subject to an inequality constraint (or a finite number thereof) of the form $g(x_1,...,x_n) \geq 0$ giving rise in general to an arbitrary (not necessarily connected) domain $\Omega^n \subset R^n$ (Fig. 1). Despite the apparent complication of the problem of evaluating sensitivity indices due to the introduction of constraints, the above formulas for main effect and total sensitivity indices remain to be applicable in any domain $\Omega^n \subset R^n$. The only modification required in Eqs. (3)-(5) is the replacement of pdf's with those defined in $\Omega^n$, which we denote in the following with a superscript '$\Omega$' (e.g. $p^\Omega(y,z)$).

For simplicity, in the following we assume that the model function $f\ x_1,...,x_n$ is defined in the unit hypercube $H^n$ (in the general case this can be achieved by applying a corresponding coordinate transformation to the original model variables) while the permissible domain $\Omega^n$ is enclosed by $H^n$: $\Omega^n \subset H^n$.

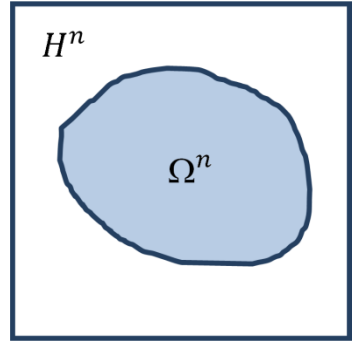

Fig. 1. Domain $\Omega^n \subset H^n$.

### 2.2. Transformation from conditional to joint and marginal distributions

We assume that there is a known procedure for the generation of vectors $(y,z)$ from a joint pdf $p^\Omega(y,z)$ defined in a non-rectangular domain $\Omega^n$, which means in particular that $p^\Omega(y,z) \equiv 0$, $(y,z) \notin \Omega^n$. Computation of $S_y$ and $S_y^T$ requires also sampling from conditional pdf's. Consider generation of a vector $(y,\bar{z})$ for the evaluation of $f(y,\bar{z})$ from a conditional probability distribution function $p^\Omega(y,\bar{z}|y)$. Using the Bayes formula $p^\Omega(y,\bar{z}|y)p^\Omega(y) = p^\Omega(y,z)$, $p^\Omega(y,\bar{z}|y)$ can be transformed as

$$p^{\Omega}(y,\overline{z}\mid y)=\frac{p^{\Omega}(y,z)}{p^{\Omega}(y)}. \tag{6}$$

Using similar transformations for $p^{\Omega}(y',\hat{z}\mid y')$ and $p^{\Omega}(\overline{y}',z\mid z)$ formulas (3)-(5) for $S_y$ and $S_y^T$ in domain $\Omega^n$ can be presented as follows

$$S_y=\frac{1}{D}\left[\int_{\Omega^s}p^{\Omega}(y)dy\left[\int_{\Omega^{n-s}}\frac{f(y,z)}{p^{\Omega}(y)}p^{\Omega}(y,z)dz\right]^2-f_0^2\right], \tag{7}$$

$$S_y=\frac{1}{D}\left[\int_{\Omega^n}f(y',z')p^{\Omega}(y',z')dy'dz'\left[\int_{\Omega^{n-s}}\frac{f(y',z)}{p^{\Omega}(y')}p^{\Omega}(y',z)dz-\int_{\Omega^n}f(y,z)p^{\Omega}(y,z)dydz\right]\right], \tag{8}$$

$$S_y^T=\frac{1}{2D}\int_{\Omega^n}\int_{\Omega^s}[f(y,z)-f(y',z)]^2\, p^{\Omega}(y,z)\frac{p^{\Omega}(y',z)}{p^{\Omega}(z)}dydy'dz. \tag{9}$$

Here the lower-dimensional integrals are understood in the sense that $\int_{\Omega^s}F(y,z)dy=\int_{(y,z)\in\Omega^n}F(y,z)dy$ and $\int_{\Omega^{n-s}}F(y,z)dz=\int_{(y,z)\in\Omega^n}F(y,z)dz$ for any $F(y,z)\in L^1(\Omega^n)$.

The usefulness of these formulas is based on the ability to sample from the joint pdf $p^{\Omega}(y,z)$ in $\Omega^n$. Although there is a technique based on the sequential sampling from an inverse cumulative distribution presented in Appendix A, it has a limited applicability. One example of an application of this technique is given in Appendix B. In the next section we present a more general approach based on the acceptance – rejection method.

### 2.3. Acceptance – rejection method

Consider the joint pdf $p(y,z)$ in the absence of constraints in $H^n$. We assume that constraining the variables to an area $\Omega^n\subset H^n$ implies that their joint pdf becomes $p^{\Omega}(y,z)$, which takes zero values in $H^n\setminus\Omega^n$ and is proportional to $p(y,z)$ within $\Omega^n$. Hence, imposing constraints does not affect the functional form of the joint pdf $p(y,z)$ inside $\Omega^n$ apart from appropriate rescaling such that $\int_{\Omega^n}p^{\Omega}(y,z)dydz=\int_{H^n}p^{\Omega}(y,z)dydz=1$. Then the 'constrained' joint pdf can be defined through the 'unconstrained' one via the relationship:

$$p^{\Omega}(y,z)=\frac{p(y,z)I^{\Omega}(y,z)}{\int_{\Omega^n}p(y,z)dydz}=\frac{p(y,z)I^{\Omega}(y,z)}{\overline{I}}, \tag{10}$$

where

$$\overline{I}=\int_{\Omega^n}p(y,z)dydz=\int_{H^n}p(y,z)I^{\Omega}(y,z)dydz \tag{11}$$

is a scaling factor and $I^{\Omega}(y,z)$ is an indicator function of the subset $\Omega^n$

$$I^{\Omega}(y,z)=\begin{cases}1, & (y,z)\in\Omega^n\\ 0, & (y,z)\notin\Omega^n.\end{cases} \tag{12}$$

The latter can be represented through the Heaviside step function $U(\cdot)$ applied to the constraint(s):

$$I^{\Omega}(y,z)=U(g(y,z)). \tag{13}$$

A constrained marginal pdf is then defined as

$$p^\Omega(y) = \int_{\Omega^n} p^\Omega(y,z)dz = \frac{1}{\overline{I}} \int_{H^{n-s}} p(y,z) I^\Omega(y,z)dz, \qquad (14)$$

while a conditional pdf takes the form

$$p^\Omega(y,\overline{z} \mid y) = \frac{p^\Omega(y,z)}{p^\Omega(y)} = \frac{p(y,z)I^\Omega(y,z)}{\int_{H^{n-s}} p(y,z)I^\Omega(y,z)dz}. \qquad (15)$$

The marginal pdf $p^\Omega(z)$ required for the evaluation of $S_y^T$ can be expressed similarly to the formula (14).

Consider now a model function $f(y,z)$ defined in $\Omega^n$. We note that it can be artificially extended to $H^n$ by defining a continuation of $f(y,z)$ such that $f(y,z) \equiv 0$ for $(y,z) \in H^n \setminus \Omega^n$ to benefit from simpler boundaries of the integration domain. Then the expected value and the variance of $f(y,z)$ are given by

$$f_0 = \int_{\Omega^n} f(y,z) p^\Omega(y,z) dy dz = \frac{1}{\overline{I}} \int_{H^n} f(y,z) p(y,z) I^\Omega(y,z) dy dz, \qquad (16)$$

$$D = \int_{\Omega^n} f^2(y,z) p^\Omega(y,z) dy dz - f_0^2 = \frac{1}{\overline{I}} \int_{H^n} f^2(y,z) p(y,z) I^\Omega(y,z) dy dz - f_0^2, \qquad (17)$$

correspondingly.

Using the expressions above for pdf's as well as the expected value and total variance of the function, the main effect and total indices can be computed through their integral formulations presented in (7)-(9). The important difference, however, is that even if the relevant pdf's are known only within the enveloping hypercube (i.e., for the unconstrained case), they can be used to directly compute their constrained counterparts on the basis of an acceptance-rejection approach invoking the indicator function of the feasible subdomain.

The DL formula (7) can be rewritten in the following simplified form

$$S_y = \frac{1}{D} \left[ \int_{H^s} \frac{\left[ \int_{H^{n-s}} f(y,z) p^\Omega(y,z) dz \right]^2}{p^\Omega(y)} dy - f_0^2 \right] \qquad (18)$$

while equations (8)-(9) remain unchanged. Note that integration in (18) is performed over lower-dimensional projections $H^s$ and $H^{n-s}$ of the unit hypercube $H^n$ as opposed to those of the domain $\Omega^n$ of potentially complex shape.

An alternative expression for the total effect index $S_y^T$ can be obtained using the identity

$$D_y^T = D - D_z \left[ E_y \, f(\overline{y}, z) \right], \qquad (19)$$

where

$$D_z \left[ E_y \, f(\overline{y}, z) \right] = \int_{\Omega^{n-s}} p^\Omega(z) dz \left[ \int_{\Omega^s} f(\overline{y}, z) p^\Omega(y, z \mid z) dy \right]^2 - f_0^2. \qquad (20)$$

Using the Bayes formula and transformations similar to presented above we obtain a DL-like formula for total indices:

$$S_y^T = 1 - \frac{1}{D} \left( \int_{H^{n-s}} \frac{\left[ \int_{H^s} f(y,z) p^\Omega(y,z) dy \right]^2}{p^\Omega(z)} dz - f_0^2 \right). \qquad (21)$$

### 2.4. Interpretation of indices

We start with Interpretation of the total index. In the case of independent inputs

$$D_y^T = D_y + D_{yz}, \quad (22)$$

where $D_{yz}$ is an interaction term between subsets $y$ and $z$. From definition (19) it follows that $D_y^T$ is the part of $D$ which remains after deduction of $D_z\left[E_y\ f(\overline{y},z)\right]$. In the case of dependent inputs $D_z\left[E_y\ f(\overline{y},z)\right]$ contains contribution to the total variance corresponding to the effect of subset $z$ on its own and its dependence with subset $y$ due to correlation or dependence via inequality constrains, which we'll call for brevity dependence contribution $D_{yz}^C$. Hence, $D_y^T$ contains dependence contribution $D_{yz}^C$ with the negative sign that is $D_y^T$ does not contain correlated/dependence contribution. On the other hand, the main effect index includes contribution corresponding to the effect of subset $y$ on its own plus dependence contribution $D_{yz}^C$. This new meaning of $D_y$ lead the authors of [6] to call $S_y$ the total correlated contribution and $S_y^T$ - the total uncorrelated contribution.

### 3. Numerical methods

In this section we present general MC estimators, acceptance-rejection estimation of Sobol' sensitivity indices using grid quadrature formulas and MC estimators based on the acceptance-rejection method.

### 3.1. MC estimators

We consider MC estimators for the evaluation of integral expressions in (7)-(9) assuming that all pdf's are defined in $\Omega^n$ and there is a sampling procedure able to generate random vectors within this domain. Formula (7) can be used to derive the DL MC estimator for $S_y$ for a single-variable index ( $y = x_i$ ). In this case $N$ points $x^{(l)}$, $l = 1, 2, ..., N$ are generated from the joint probability distribution $p^\Omega(y,z)$. For each random variable $y = x_i$, the sample set $x^{(l)}$, $l = 1, 2, ..., N$ is sorted in the ascending order on the interval $[0,1]$ with respect to this variable and then subdivided into $N_y$ equally populated partitions (bins) each containing $N_z = N/N_y$ points ($N_y < N$). Within each bin we calculate the local mean value $E_z\left[f(y,z)|y_j^A\right] \approx \dfrac{1}{N_z}\sum_{k=1}^{N_z} f(y_{j_k}, z_{j_k})$, where $j$ is the index of a bin containing $N_z$ points, $j = 1, ..., N_y$ and $y_j^A$ is a mean value of $y$ within this bin. We note, that $y_j^A$ is not actually computed and it is used only for notation purposes. Finally, the variance of all conditional averages is computed to yield the following double loop reordering (DLR) formula:

$$S_y \approx \frac{1}{D}\left[\frac{1}{N_y}\sum_{j=1}^{N_y}\frac{\left(\frac{1}{N_z}\sum_{k=1}^{N_z}f(y_{j_k},z_{j_k})\right)^2}{\frac{1}{N_z}\sum_{k=1}^{N_z}p^{\Omega}(y_{j_k})}-f_0^2\right] = \frac{1}{D}\left[\frac{1}{N}\sum_{j=1}^{N_y}\frac{\left(\sum_{k=1}^{N_z}f(y_{j_k},z_{j_k})\right)^2}{\sum_{k=1}^{N_z}p^{\Omega}(y_{j_k})}-f_0^2\right]. \quad (23)$$

The subdivision into bins is done in the same way for all inputs using the same set of sample points. A critical issue is the link between $N$ and $N_y$. It was suggested in [16] to use as a "rule of thumb" $N_y \approx \sqrt{N}$.

We note that it is practically impossible to generalise the DLR formula for more than one index, hence there is no similarly efficient DLR formula which allows to compute total Sobol' sensitivity indices.

The expected value and total variance are computed using the MC estimators

$$f_0 \approx \frac{1}{N}\sum_{l=1}^{N}f(y_l,z_l). \quad (24)$$

$$D \approx \frac{1}{N}\sum_{l=1}^{N}f^2(y_l,z_l)-f_0^2. \quad (25)$$

The MC estimator for the main effect index based on formula (8) has a form

$$S_y \approx \frac{1}{DN}\sum_{l=1}^{N}\left(f\left(y_l',z_l'\right)\left(\frac{f(y_l',z_l)}{p^{\Omega}(y_l')}-f\left(y_l,z_l\right)\right)\right), \quad (26)$$

and an estimator for the total effect index based on formula (9) can be written as:

$$S_y^T \approx \frac{1}{2DN}\sum_{l=1}^{N}\frac{1}{p^{\Omega}(z_l)}\left(f\left(y_l,z_l\right)-f\left(y_l',z_l\right)\right)^2. \quad (27)$$

The usefulness of these estimators relies on the ability to sample from the joint pdf $p^{\Omega}(y,z)$ in $\Omega^n$. In the next two sub-sections we present a practical approach for such sampling based on the acceptance-rejection method.

### 3.2. Acceptance-rejection estimation of sensitivity indices using grid quadrature formulas

Integrals in (16)-(18) and (21) can be efficiently evaluated using grid quadrature methods for problems of low and medium dimensionality. Application of grid quadrature methods is rather straightforward when integration domains are hyperrectangular as in the case of the acceptance-rejection approach presented in section 2.3. However, given that the feasible domain $\Omega^n$ can be arbitrary in shape, efficient higher-order quadrature methods (such as those based on Gaussian or Clenshaw–Curtis quadrature [17]) lose their advantages as the integrands are not differentiable at the boundary of the permissible domain $\Omega^n$ (the indicator function has a jump discontinuity at the boundary of $\Omega^n$). Thus the use of lower-order integration formulas such as the second-order multidimensional trapezoidal rule used in this paper is preferable because they are less sensitive to non-smooth or discontinuous integrands.

Numerical integration can gain additional efficiency by performing preliminary bracketing of the domain $\Omega^n$ (or finding its 'minimum bounding box') within $H^n$ to minimize the number of rejected points during sampling. This is especially important in higher dimensions, when the number of grid points in each dimension can't be large due to the "curse of dimensionality". The total number of grid points is $N=k^n$, where $k$ is the number of points in each dimension. $k$ includes both internal and boundary points and it can't be smaller than 4. We don't consider the case of the different values of k for different directions due the lack of prior information about importance of individual inputs in the general case.

The bracketing is performed by first using a uniform grid in $H^n$ to evaluate new lower ($x_i^{\min} \geq 0$) and upper ($x_i^{\max} \leq 1$) bounds for each variable. $x_i^{\min}$ and $x_i^{\max}$ are defined as

$$x_i^{\min} = \inf \, x_i : \int_{H^{n-1}} I^{\Omega}(x) dx_{-i} > 0 \,, \tag{28}$$

$$x_i^{\max} = \sup \, x_i : \int_{H^{n-1}} I^{\Omega}(x) dx_{-i} > 0 \,, \tag{29}$$

where $x_{-i}$ is the vector of all model variables but $x_i$: $x_{-i} = (x_1,...,x_{i-1},x_{i+1},...,x_n)$. Once the new bounds have been determined the found hyperrectangle $H_{\min}^n = \prod_{i=1}^{n} \left[ x_i^{\min}, x_i^{\max} \right]$ which can be seen as a tight enclosure of $\Omega^n$ is covered with a new grid, which is used to sample values of the model function $f(x)$ at points where the indicator function $I^{\Omega}(x)$ is non-zero.

This approach was successfully implemented for problems of dimensionality up to 10. The results are presented in Section 4. The use of this approach for higher dimensional models is computationally prohibitive. This is also the main reason why the grid quadrature approach is not applicable to the modified formulas (8)-(9) since the effective number of dimensions for integration in this case is $(3n-s)$ and $(n+s)$, respectively. Extension to higher dimensions may be possible with application of sparse grid methods.

### 3.3. MC estimators for the acceptance-rejection method

In this subsection we consider MC estimators of $S_y$ and $S_y^T$ valid when $p(y,z)$ is known but the constrained pdf $p^{\Omega}(y,z)$ is not known explicitly. We also assume that there is no explicit technique for sampling in $\Omega^n$ of arbitrary shape. Firstly we derive an MC estimator for $S_y$ defined by (18). Its approximation leads to the DLR approach discussed in subsection 3.1. Based on a sequence of points $(y_l, z_l)$, $l = 1,...,N$ sampled from the joint pdf $p(y,z)$, the mean value $f_0$ and total variance $D$ of $f(y,z)$ are estimated as

$$f_0 \approx \frac{1}{\bar{I} N} \sum_{l=1}^{N} f(y_l, z_l) I^{\Omega}(y_l, z_l), \tag{30}$$

$$D \approx \frac{1}{\bar{I} N} \sum_{l=1}^{N} f(y_l, z_l) - f_0^{\,2} \, I^{\Omega}(y_l, z_l), \tag{31}$$

where

$$\bar{I} \approx \frac{1}{N} \sum_{l=1}^{N} I^{\Omega}(y_l, z_l) \tag{32}$$

is the scaling factor.

Further we assume that $y = x_i$. The integral $F(y) = \int_{H^{n-s}} f(y,z) p^{\Omega}(y,z) dz$ in the numerator of the integrand of (18) can be approximated by the following estimator

$$F(y_j^A) \approx \frac{1}{\bar{I} N_z} \sum_{k=1}^{N_z} f(y_{j_k}, z_{j_k}) I^\Omega(y_{j_k}, z_{j_k}), \tag{33}$$

i.e., by averaging within each of $N_y$ bins, where $j$ is the index of a bin and $y_j^A$ represents an average value of $y$ in the $j$-th bin. The marginal distribution function $p^\Omega(y)$ in the denominator of (18) is approximated within each bin by

$$p^\Omega(y_j^A) \approx \frac{1}{\bar{I} N_z} \sum_{k=1}^{N_z} I^\Omega(y_{j_k}, z_{j_k}). \tag{34}$$

Combining the above we obtain a DLR MC estimator for the main effect index defined by formula (18) in the following form:

$$S_y \approx \frac{1}{D} \left[ \frac{1}{N_y} \sum_{j=1}^{N_y} \frac{F^2(y_j^A)}{p^\Omega(y_j^A)} - f_0^2 \right]. \tag{35}$$

An MC estimator for $S_y$ based on the modified formula (8) has the form:

$$S_y = \frac{1}{\bar{I}^2 DN} \sum_{l=1}^{N} \left( f(y_l', z_l') I(y_l', z_l') \left( \frac{f(y_l', z_l) I(y_l', z_l)}{p^\Omega(y_l')} - f(y_l, z_l) I(y_l, z_l) \right) \right). \tag{36}$$

Similarly an estimator for the total effect index (9) can be written as:

$$S_y^T = \frac{1}{2\bar{I}^2 DN} \sum_{l=1}^{N} \left( f(y_l, z_l) I(y_l, z_l) - f(y_l', z_l) I(y_l', z_l) \right)^2 \frac{1}{p^\Omega(z_l)}, \tag{37}$$

where the estimation of the marginal probability distribution $p^\Omega(z)$ requires additional sampling as described below.

The marginal distribution $p^\Omega(y)$ approximated by (38) is estimated by subdividing the whole set of sample points $(y_l, z_l)_{l=1}^{N}$ into $N_y$ bins according to the values of $y_l$ so that the probability of points ending up in any of the bins is the same and equal to $\Delta p_y = 1/N_y$ (we recall that $y = x_i$ is scalar, $y \in [0,1]$). To keep the same resolution in the estimation of $p^\Omega(z)$ an $(n-s)$-dimensional "$z$-bin" can be defined by choosing probabilities $\Delta p_{z_i} = 1/N_z$ for each $i = 1, \ldots, n-s$. Owing to the "rule of thumb" introduced earlier $N_z \approx N_y \approx \sqrt{N}$ thus leading to similar probabilities defining a bin in all dimensions, whether $y$ or $z_i$, $i = 1, \ldots, n-s$. However, such a definition results in a dramatically decreasing number of points of the original sample falling within such a bin when the model dimensionality $n$ (and hence that of the bin, $n-s$) increases. Indeed, the probability that a random point sampled from the joint pdf falls within a $z$-bin defined in such a way is given by $\Delta p_{\text{bin},z} = \prod_{i=1}^{n-s} \Delta p_{z_i} = N_z^{-(n-s)} = N^{-\frac{n-s}{2}}$ so that the number of points of the original sample expected to be found in such a bin is $N_{\text{bin},z} = N \Delta p_{\text{bin},z} = N^{-\frac{n-s}{2}+1}$, which is greater than one only if $n-s < 2$ or, in the most typical case of $s=1$ considered here, if $n < 3$. Therefore, to estimate $p^\Omega(z_l)$ with acceptable accuracy a relatively small additional number of samples of the form $(y_q, z_l)$, $q = 1, \ldots, N_z'$ is required which yields the estimator

$$p^\Omega(z_l) \approx \frac{1}{\bar{I} N_z'} \sum_{q=1}^{N_z'} I^\Omega(y_q, z_l). \tag{38}$$

We note that only the indicator function is evaluated at each of these points and although this additional sampling at each $z_l$ increases the computational cost of the estimator (37), this increase may not be detrimental if the computational cost of the evaluation of model constraints is significantly smaller than that of the model function itself (which is often the case in practical problems). A value $N'_z = 64$ was used in the test cases reported below, which was found to be high enough to provide good convergence rates.

There are two different ways of computing the set $\{S_i, S_i^T\}$, $i=1,...,n$. The first one is based on the modified formulas (36), (37). The required number of function evaluations in this case is $N_{CPU} = N(n+2)$. The second one is based on using the DLR formula (35) for computation $S_i$ (which requires $N_{CPU} = N$ function evaluations) and formula (37) for computation $S_i^T$ (which requires $N_{CPU} = N(n+2)$ function evaluations). However, when (35) is used in conjunction with the modified formula (37) for total effect sensitivity indices, the former can benefit from using all the available $N_{CPU} = N(n+2)$ samples required to compute $S_i^T$ to enhance the accuracy of $S_i$ estimation. All derived MC estimators can be computed using MC or QMC sampling.

## 4. Test cases

For all test cases considered in this section we found analytical values of sensitivity indices, so that they can be used as benchmarks for verification of numerical estimates.

### 4.1. Product function in triangle domain $\Omega_1$.

Consider the function
$$f(x_1, x_2) = x_1 x_2 \tag{39}$$
with both variables uniformly distributed in an upper triangle $\Omega_1$ shown in Fig. 2 with joint pdf
$$p^{\Omega_1}(x_1, x_2) = \begin{cases} 2, & (x_1, x_2) \in \Omega_1 \\ 0, & (x_1, x_2) \notin \Omega_1 \end{cases}. \tag{40}$$
$\Omega_1$ is defined by constraint $g(x_1, x_2) = x_1 + x_2 - 1 \geq 0$.

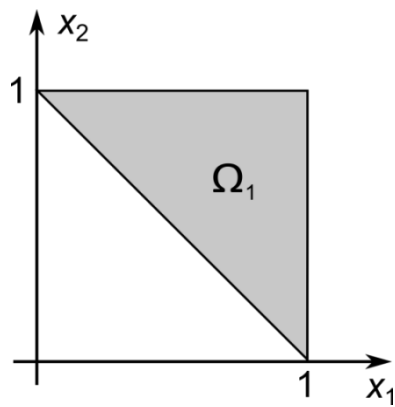

Fig. 2. Definition of area $\Omega_1$.

Using the Bayes theorem the marginal pdf of $x_1$ and the pdf of $x_2$ conditional on $x_1$ can be computed explicitly:

$$p_1^{\Omega_1}(x_1) = \int_{1-x_1}^{1} p^{\Omega_1}(x_1, x_2)\, dx_2 = \int_{1-x_1}^{1} 2\, dx_2 = 2x_1, \tag{41}$$

$$p_2^{\Omega_1}(x_2 \mid x_1) = \frac{p^{\Omega_1}(x_1, x_2)}{p_1^{\Omega_1}(x_1)} = \frac{1}{x_1}. \tag{42}$$

Expectation and total variance for this function can also be computed explicitly:

$$f_0 = E\{f(x_1, x_2)\} = \int_0^1 \int_{1-x_1}^{1} x_1 x_2\, 2\, dx_2\, dx_1 = \frac{5}{12}, \tag{43}$$

$$D = \int_0^1 \int_{1-x_1}^{1} (x_1 x_2)^2\, 2\, dx_2\, dx_1 - \left(\frac{5}{12}\right)^2 = \frac{3}{80}. \tag{44}$$

Using the definition (7) of the main effect index in a non-rectangular area and the relationship between the main effect and total indices $S_{x_2}^T = 1 - S_{x_1}$ the sensitivity indices for the product function (39) are evaluated as follows:

$$S_{x_1} = \frac{1}{D}\left[\int_0^1 p_1^{\Omega_1}(x_1)\, dx_1 \left(\int_{1-x_1}^{1} f(x_1, x_2) p_2^{\Omega_1}(x_2 \mid x_1)\, dx_2\right)^2 - f_0^2\right] =$$

$$= \frac{80}{3}\left[\int_0^1 2x_1\, dx_1 \left(\int_{1-x_1}^{1} x_1 x_2 \frac{1}{x_1}\, dx_2\right)^2 - \left(\frac{5}{12}\right)^2\right] = \frac{7}{27}, \tag{45}$$

$$S_{x_2}^T = 1 - S_{x_1} = \frac{20}{27}. \tag{46}$$

Owing to the symmetry of the function and the area $\Omega_1$ $S_{x_2} = S_{x_1}$ and $S_{x_1}^T = S_{x_2}^T$.

Note that in the absence of the constraint, i.e. when the variables are uniformly distributed in the unit square, the sensitivity indices are $S_{x_1} = S_{x_2} = 3/7$ and $S_{x_1}^T = S_{x_2}^T = 4/7$, while the mean value and total variance are $f_0 = 1/4$ and $D = 7/144$.

### 4.2. *g*-function in triangular domains $\Omega_1$ and $\Omega_2$.

Consider the so-called *g*-function which is often used in GSA for illustration purposes [1]:

$$f = \prod_{i=1}^{2} \frac{|4x_i - 2| + a_i}{1 + a_i} \tag{47}$$

with parameters $a_1 = 0$ and $a_2 = 1$ (see Fig. 3) and $x_1, x_2$ uniformly distributed in the unit square. Below we compute its sensitivity indices when the feasible domain is restricted by a parametric linear (Fig. 4) or nonlinear (Fig. 9) constraint.

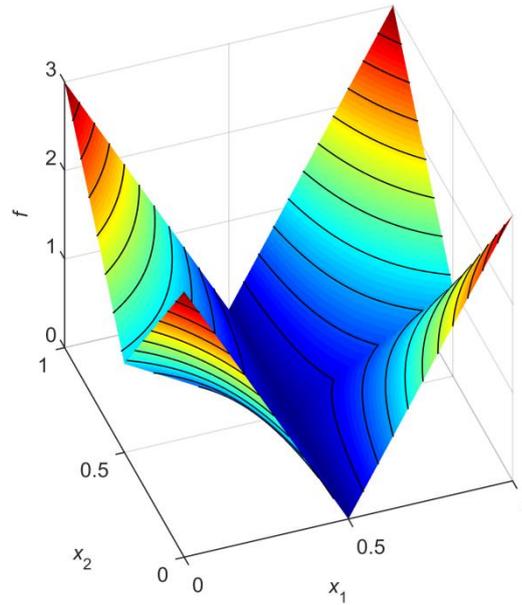

Fig. 3. Plot of the *g*-function (47)

### 4.2.1. A linear constraint

In this subsection we consider a more general parametric linear constraint (illustrated in Fig. 4) than that used in the test case presented in the previous subsection. The shape and size of the feasible domain $\Omega_2$ are defined by the variable angle $\alpha$ between the top side of the unit square and the line defined by the linear constraint

$$g(x_1, x_2) = 1 - \tan(\alpha) x_1 - x_2 \geq 0 .\qquad(48)$$

Firstly we obtain analytical results for two values of $\alpha$: $\pi/6$ and $\pi/4$. While the latter value leads to a symmetrical domain (complementary to $\Omega_1$ in Fig. 2), the former does not result in any specific simplification of the problem. The reference values of the mean, total variance and main effects are given in Table 1. Exact solution for $\alpha = \pi/6$ was obtained using symbolic integration in Maple®, however the resulting expressions are too cumbersome to report them here, so the values are given in the decimal form to sufficient precision for error analysis (see below). As noted above, the values of total sensitivity indices are readily obtained from the relationships $S^T_{x_2} = 1 - S_{x_1}$, $S^T_{x_1} = 1 - S_{x_2}$ which are valid in the 2D case.

Figs. 5, 6 show the variations of the expected value, total variance and the main effect and total sensitivity indices with the angle $\alpha$ changing from 0 to $\pi/2$, which corresponds to the whole range from the completely unconstrained 2D problem to the degenerate 1D case $x_1 = 0$, $x_2 \in \langle 0,1 \rangle$. Analytical exact values given in Table 1 are denoted by overlaid symbols. The results perfectly respect limiting cases given in Table 1. The numerical results were obtained with the use of the grid quadrature method.

**Table 1**. Exact values of mean, total variance and main effects for selected values of $\alpha$

| $\alpha$ | $f_0$ | $D$ | $S_{x_1}$ | $S_{x_2}$ |
|---|---|---|---|---|
| $\pi/6$ | 0.9714128589 | 0.4483218079 | 0.7703487112 | 0.3214987100 |
| $\pi/4$ | 1 | 4/9 | $-\dfrac{93}{40} + \dfrac{9}{2} \ln 2$ | $-\dfrac{9}{20} + \dfrac{9}{8} \ln 2$ |

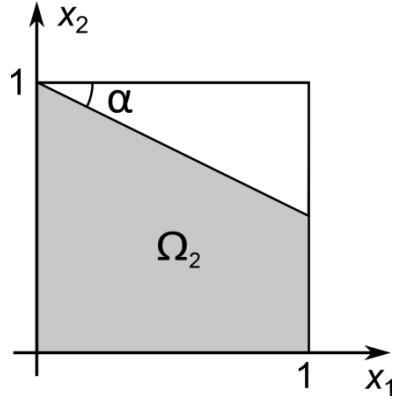

Fig. 4. Domain $\Omega_2$ as defined by the linear constraint (48).

As the value of $\alpha$ increases the levels of importance of the variables swap with $x_1$ being initially significantly more important. However, as $\alpha \to \pi/2$ the domain $\Omega_2$ degenerates into the segment $x_1 = 0$, $0 \le x_2 \le 1$ so that the model function clearly ceases to depend on $x_1$. This is reflected by the value of the total effect index $S^T_{x_1}$ dropping to zero and a simultaneous increase of $S^T_{x_2}$ towards unity.

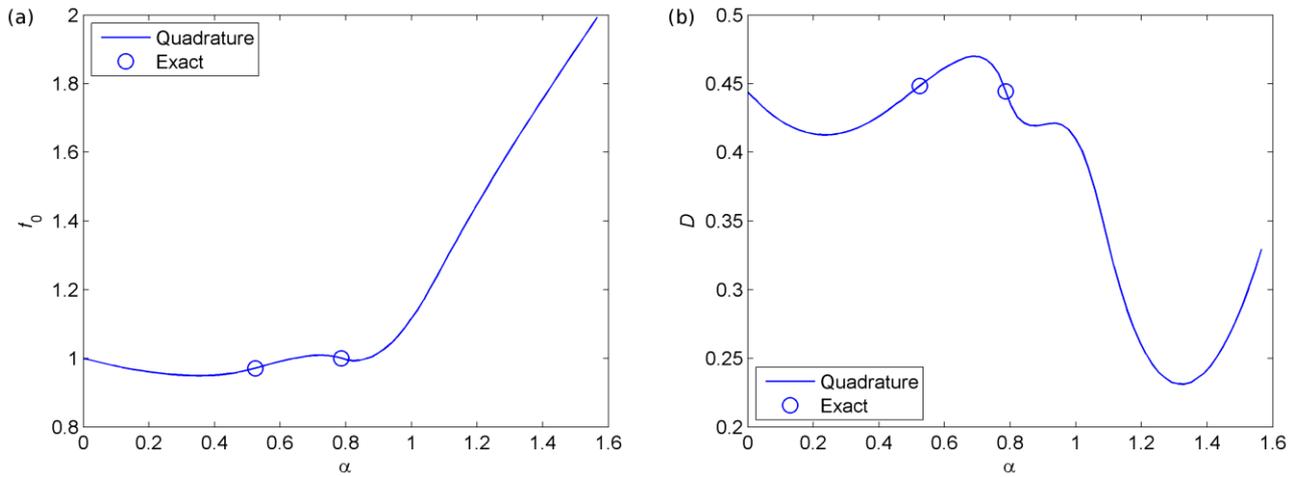

Fig. 5. Values of (a) $f_0$ and (b) $D$ for the 2D $g$-function (47) versus angle $\alpha$.

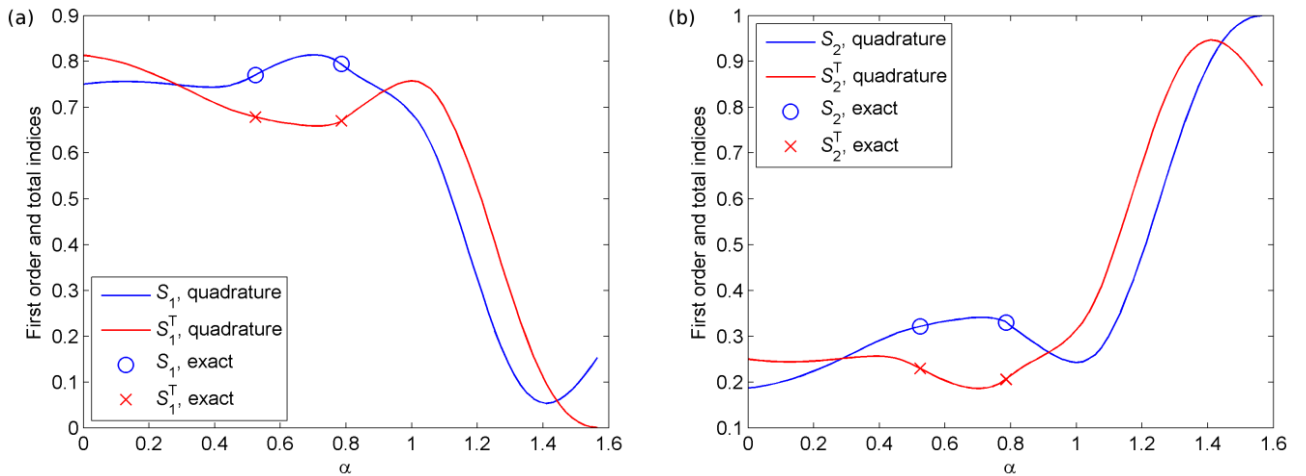

Fig. 6. Values of (a) $S_{x_1}$ and $S_{x_1}^T$, and (b) $S_{x_2}$ and $S_{x_2}^T$ for the 2D g-function (47) versus angle $\alpha$.

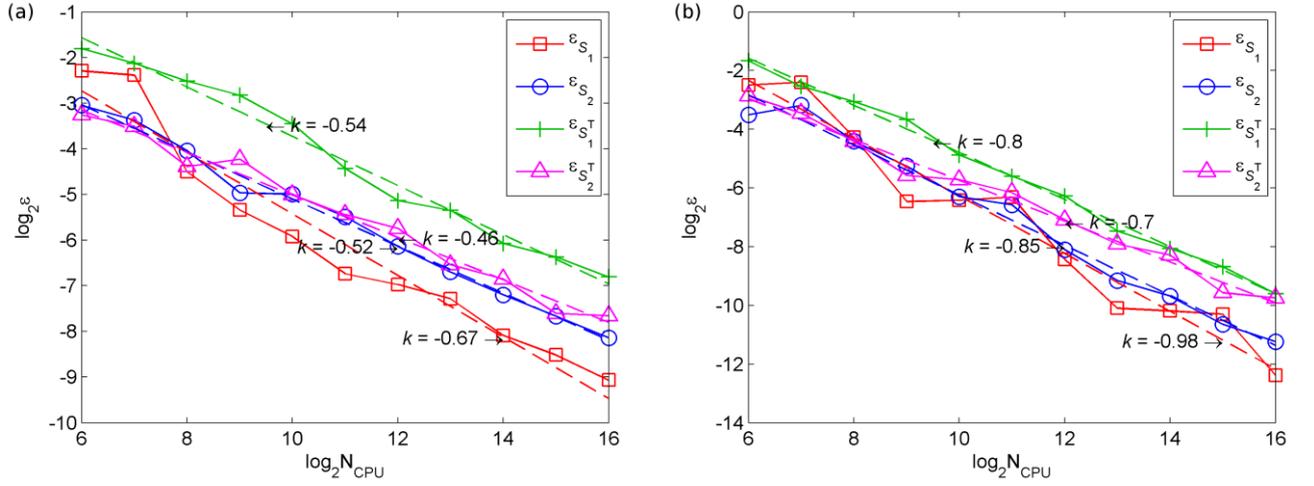

Fig. 7. Convergence of the MC estimators of $S_i$ (35) and $S_i^T$ (37) for the 2D g-function (47) defined in $\Omega_2$ for $\alpha = \pi/6$ obtained using (a) MC (b) QMC methods.

We also tested the performance of the MC estimators. The root-mean-square error (RMSE) obtained using the MC estimator (35) for the main effect indices and the DLR estimator (37) for total effects indices for a fixed value of $\alpha = \pi/6$ is presented in Fig. 7. To reduce the scatter in the error estimation the values of RMSE were averaged over L = 50 independent runs:

$$\varepsilon_i = \left( \frac{1}{L} \sum_{l=1}^{L} \left( \frac{I_{i,l}^* - I_0}{I_0} \right)^2 \right)^{\frac{1}{2}}.$$

Here $I_i^*$ is the numerical value of a particular estimator, $I_0$ is the corresponding analytical value. The RMSE is approximated by a trend line $cN^k$. Values of $k$ are given on the plots. The convergence rates of $S_i$ and $S_i^T$ for the MC method (Fig. 7a) are close to theoretically expected value of 1/2. On the other hand, the convergence rates for the QMC method (Fig. 7b) are significantly higher (not worse than 0.85 for the $S_i$ and 0.7 for $S_i^T$). It should be noted that both approaches suffer from decreased convergence rates when the area of $\Omega_2$ diminishes (data not shown) and the domain becomes elongated along the $x_2$ axis. However, even in this case the QMC method still outperforms MC.

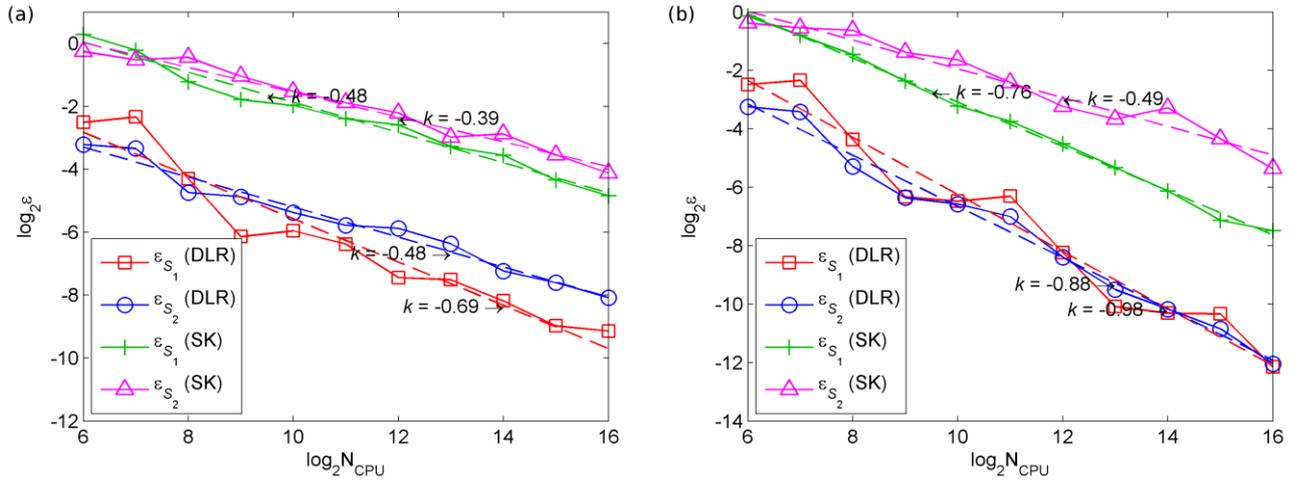

Fig. 8. Convergence of the DLR (35) and SK (36) MC estimators for $S_i$ for the 2D *g*-function (47) defined in $\Omega_2$ for $\alpha = \pi/6$ using (a) MC (b) QMC methods.

Finally, we compare the performance of the two presented estimators for the main effect indices: namely, the DLR (35) and the estimator (36) for the modified formula (8) by S. Kucherenko et al [11], the latter is denoted as SK. The results presented in Fig. 8 show significant advantages of the DLR approach when considering the convergence rate versus the total number of sampling points $N_{\text{CPU}}$.

### 4.2.2. A parabolic constraint

In this subsection we consider a more complex example involving a nonlinear constraint of the form:

$$g(x_1, x_2) = x_2 - \beta x_1 (1 - x_1) \geq 0, \qquad (49)$$

which describes the part of the unit square above a parabola (Fig. 9). Depending on the value of $\beta$ the permissible domain is either connected ($\beta \leq 4$) or disconnected ($\beta > 4$) as illustrated in Fig. 9. We note that the domain is non-convex for any value of $\beta$.

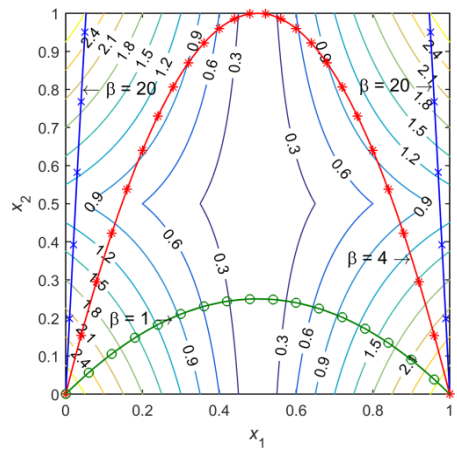

Fig. 9. A series of parabolic constraints defined by (49) for different values of $\beta$ and the contour lines of the *g*-function (47). See Fig. 3 for a 3D plot of the *g*-function

Variations of the function mean, total variance and main and total effect sensitivity indices reveal highly nonlinear behaviour shown in Figs. 10 and 11. The dashed lines in these Figs indicate the critical value $\beta = 4$ for the connectedness of the feasible domain. It is also worth noting that the limits $\beta \to 0$ and $\beta \to \infty$ represent the two extreme cases corresponding either to the unconstrained situation (unit square) or to the degenerate case when the feasible domain is the union of two disjoint segments $x_1 = 0, 1$, $x_2 \in 0,1$, respectively. In the latter limit the numerical values of $f_0$, $D$, main and total effect sensitivity indices tend toward the same limit as those under a linear constraint when $\alpha \to \pi/2$ (compare with Figs. 5, 6).

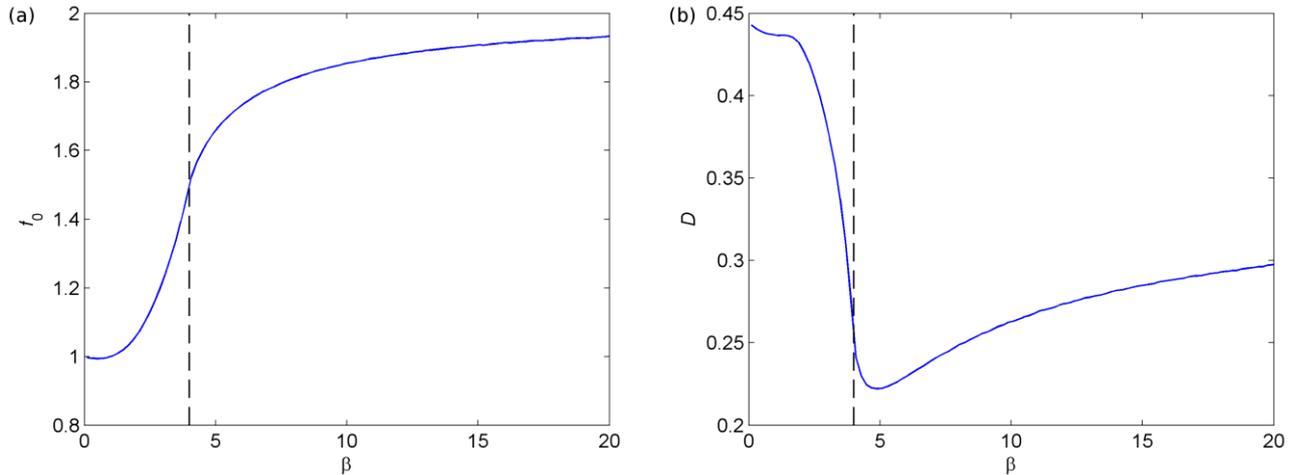

Fig. 10. Values of (a) $f_0$ and (b) $D$ for the 2D g-function (47) versus parameter $\beta$ from parabolic constraint (49).

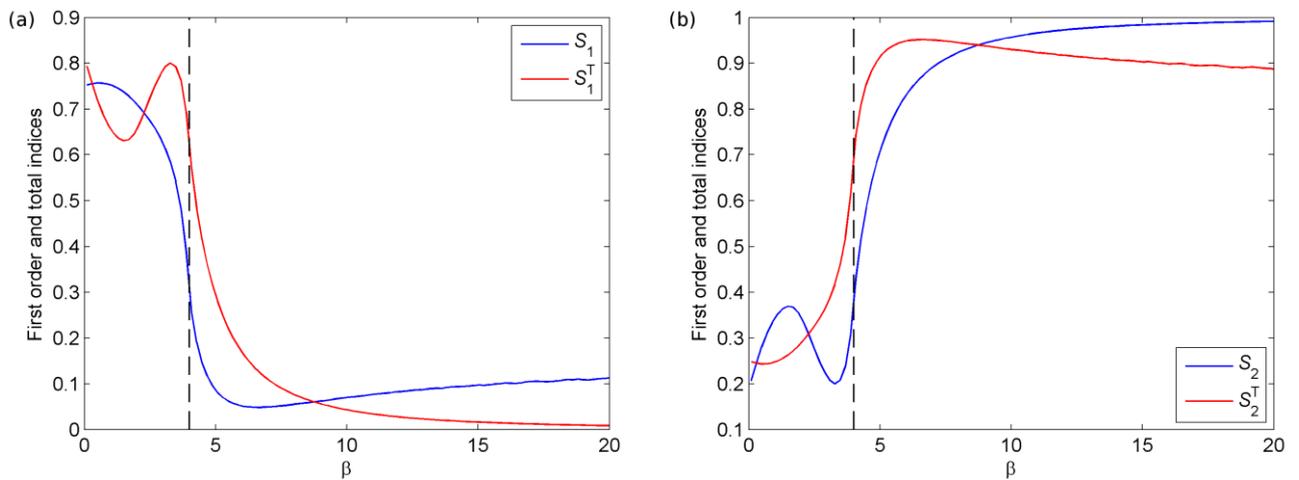

Fig. 11. Values of (a) $S_{x_1}$ and $S_{x_1}^T$, and (b) $S_{x_2}$ and $S_{x_2}^T$ for the 2D g-function (47) versus parameter $\beta$ from parabolic constraint (49).

Consider the particular case of $\beta = 4$ which is the smallest value of $\beta$ when the feasible domain is essentially disconnected. The exact values for the mean, total variance and main effect indices are:

$$f_0 = 3/2,$$

$$D = \frac{521}{1260} - \frac{4}{35}\sqrt{2},$$

$$S_{x_1} = \frac{3}{144\sqrt{2}-521},$$

$$S_{x_2} = -\frac{3}{2}\frac{264\sqrt{2}-293}{144\sqrt{2}-521}. \tag{50}$$

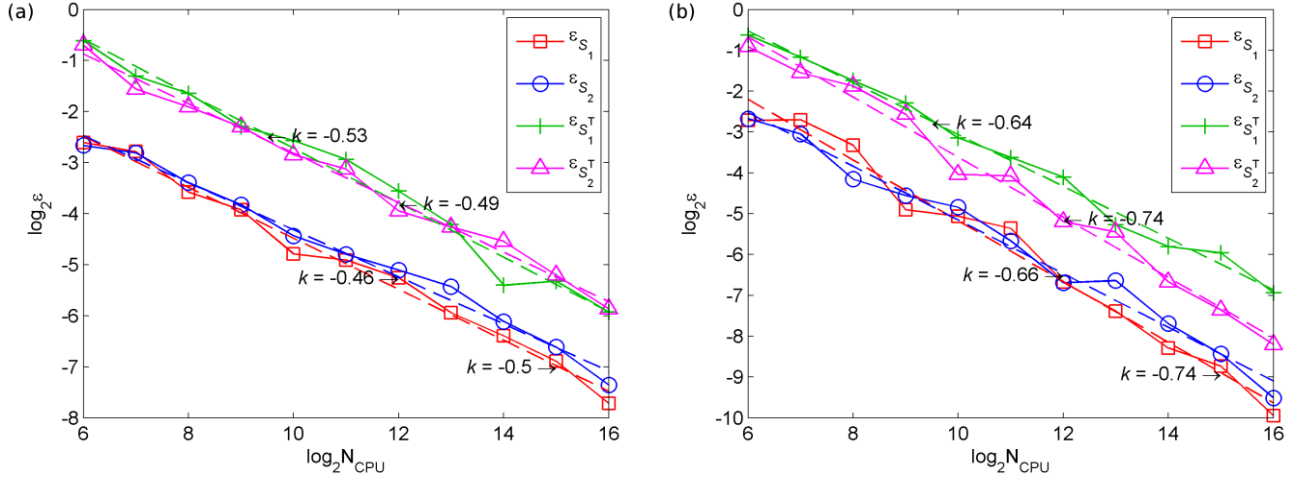

Fig. 12. Convergence of the MC estimates of $S_i$ (35) and $S_i^T$ (37) for the 2D *g*-function (47) under the parabolic constraint (49) with $\beta = 4$ obtained using (a) MC (b) QMC methods.

Fig. 12 shows the variation of RMSE versus the total number of sampling points $N_{\text{CPU}}$ for MC and QMC methods. The latter is clearly superior to MC.

### 4.3. K-function

K-function is defined as

$$K = \sum_{i=1}^{n}(-1)^i \prod_{j=1}^{i} x_j, \tag{51}$$

where variables $x_j$, $j=1,...,n$ are independent uniformly distributed random variables in [0, 1]. K-function is also used in GSA for illustration purposes (see f.e.[ 18]).

We consider four different cases for domain definitions. The first one is an unconstrained problem ( $x \in H^n$ ). In the other three cases the unit hypercube is divided by a hyperplane into two parts one of which is the permissible region for the problem variables $x_j$, $j=1,...,n$. All of these cases are considered in the four-dimensional space $n = 4$. The constraints are as follows:

$$I_1 : x_1 + x_2 \leq 1, \tag{52}$$

$$I_2 : x_3 + x_4 \leq 1, \tag{53}$$

$$I_3 : x_1 + x_3 \leq 1. \tag{54}$$

These constraints can be represented using the following indicator functions: $I_1 = U(1-x_1-x_2)$, $I_2 = U(1-x_3-x_4)$, $I_3 = U(1-x_1-x_3)$. See Fig. 13 for a schematic plot illustrating $I_1$ constraint in the 3D space.

For numerical estimates we make use of grid quadrature formulas (multidimensional trapezoidal rule) presented in subsection 3.2. In order to assess the accuracy of numerical computations in the unconstrained case the exact solution for the total effect indices reported in [19] was used while the analytical solution for the main effect sensitivity indices was derived in this paper:

$$S_i = \frac{\left(\frac{1}{2}\right)^{2i-2} + \left(-\frac{1}{2}\right)^{n+i-2} + \left(\frac{1}{2}\right)^{2n}}{\frac{3}{2} - \frac{3}{5}(-1)^n \left(\frac{1}{2}\right)^{n-1} + \frac{1}{10}\left(\frac{1}{3}\right)^{n-3} - 3\left(\frac{1}{2}\right)^{2n}}. \tag{55}$$

For the situations involving the constraints (52)-(54) exact values of $S_i$ and $S_i^T$ were obtained with the aid of Maple® software.

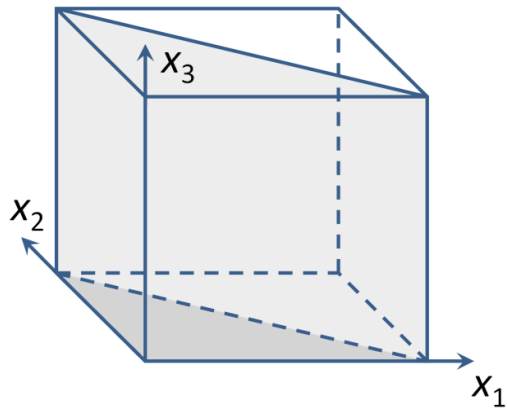

Fig. 13. Schematic representation of permissible regions for the K-function (shaded area) in the 3D case.

The exact values of $S_i$ and $S_i^T$ for all four cases are presented in Fig. 14. For the unconstrained case the most influential input is $x_1$ followed by $x_2$ with $x_3$ and $x_4$ having equal and much less significant contributions to the function variance.

Compared to the unconstrained case, the introduction of the first constraint defined by (52) (the indicator function $I_1$) leads to a significant increase of the main effect indices $S_1$ and $S_2$ accompanied by a simultaneous decrease of $S_3$ and $S_4$. It reflects the fact that when the constraint (52) is imposed an even more substantial part of the variance is contributed by the first two terms of the K-function compared to the unconstrained case.

The second constraint (53) (the indicator function $I_2$) introduces additional interaction between the less important variables $x_3$ and $x_4$. Since $x_1$ and $x_2$ are not affected by $I_2$ their sensitivity indices (both main and total effects) remain at the same level as in the unconstrained case. However, variable $x_3$ features more prominently than $x_4$ compared with the unconstrained case.

The last constraint (54) (the indicator function $I_3$) links an influential variable $x_1$ and a significantly less important $x_3$. It results in the decreased values of sensitivity indices of $x_1$, while the values of sensitivity indices of $x_2$ are increased. On the other hand, the main effect index of $x_3$ becomes significantly more important than in the unconstrained case drastically outrunning $x_4$. However, the corresponding total effects $S_3^T$ and $S_4^T$ turn out to be the same.

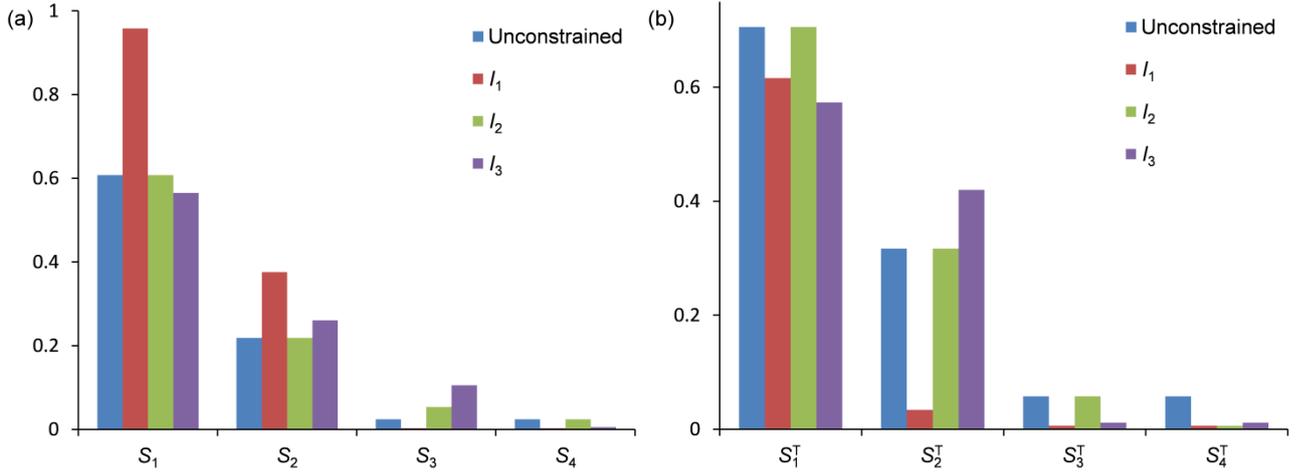

Fig. 14. (a) Main effect and (b) total sensitivity indices of the K-function in $H^4$ for the unconstrained and constraints (52)-(54) cases.

Fig. 15 shows the performance of the grid quadrature integration method. Absolute errors of the estimates of $S_i$ (Fig. 15,a) and $S_i^T$ (Fig. 15,b) decrease in all cases at a rate $\varepsilon \sim O(N^{-1/2})$. This is indeed as expected since the error of the trapezoidal rule (regardless of the dimensionality) is of the form $\varepsilon \sim O(h^2)$, where $h$ is the grid spacing assumed to be equal in all dimensions. Since the number of function evaluations in an $n$-dimensional space $N \sim (1+1/h)^n$ (with equality instead of proportionality in the unconstrained case) we obtain that in general $\varepsilon \sim O(N^{-2/n})$.

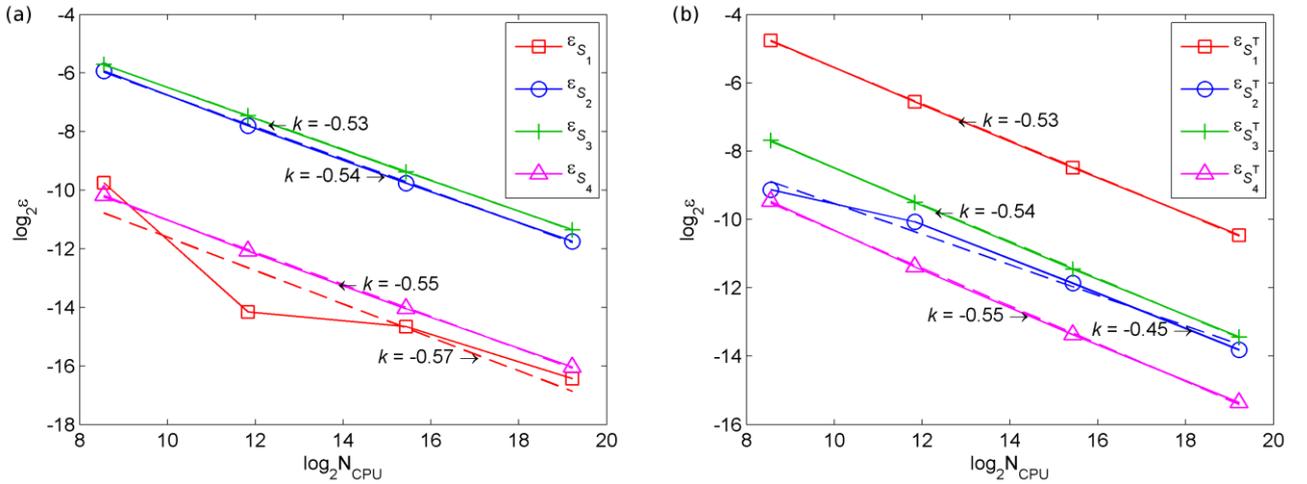

Fig. 15. Relative error in (a) $S_i$ and in (b) $S_i^T$ evaluated using grid quadrature versus the number of function evaluations $N_{CPU}$, for the K-function under constraint (54) (defined by the indicator $I_3$).

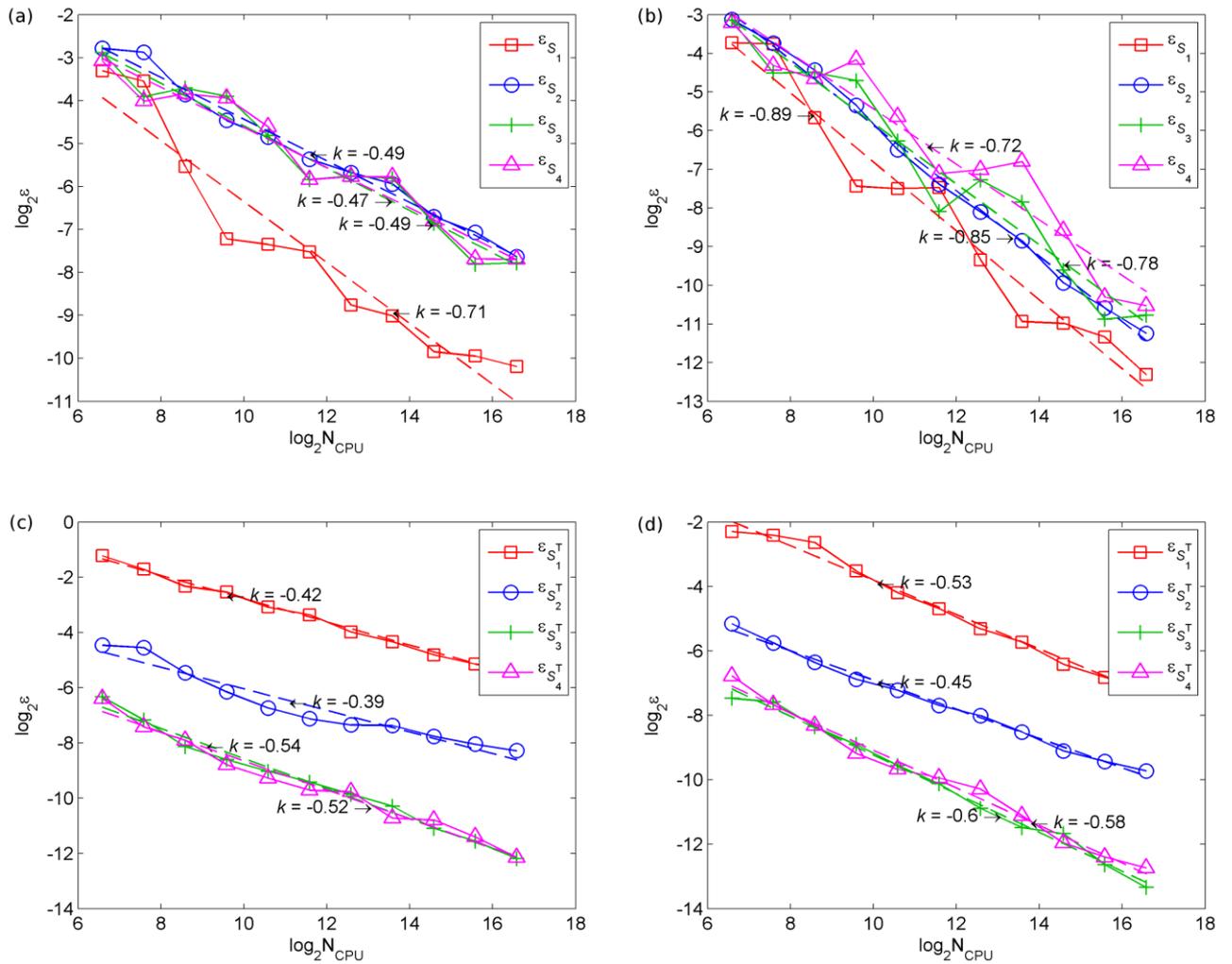

Fig. 16. Convergence rates of the (a, b) DLR estimator (35) of $S_i$ and (c, d) modified Sobol' estimator (37) of $S_i^T$ for the 4D K-function under constraint (52) (defined by the indicator $I_1$) obtained using (a, c) MC and (b, d) QMC.

Fig. 16 illustrates the performance of the MC estimators for DLR approach for the main effect and the modified Sobol' formula for the total effect indices. Similarly to the previous examples, the use of QMC sampling results in faster convergence.

**Conclusions**

In this work, we have proposed a novel concept of constrained GSA which adds the ability to analyse model output variance in arbitrarily shaped *n*-dimensional domains. This amounts to greatly expanding the scope of GSA by allowing model variables to be subject to inequality constraints, which is common in a range of situations of practical importance.

The proposed formulas build upon Sobol' sensitivity indices and their recent development for models with dependent variables [11]. The advantage of the presented formulations is that no prior knowledge of conditional or marginal distributions is assumed. All the required dependences are derived from the joint pdf in the presence of constraints. It is shown that the knowledge of the joint pdf corresponding to the unconstrained formulation is sufficient to build numerical estimators for sensitivity indices.

Three types of numerical estimators of the sensitivity indices have been proposed. Grid quadrature may compete with MC estimators for low and medium dimensional models. However, its convergence rate rapidly degrades with increasing model dimensionality. Despite its simple concept, the DLR approach demonstrates good convergence when applied to the evaluation of main effects sensitivity indices outperforming the MC estimator of the modified Sobol' formulas. On the other hand, DLR is not a viable approach for the evaluation of total effect indices for which the modified Sobol' formula gives good results.

Further work is needed to develop a clear interpretation of the cGSA results, in particular by decomposing the variance contribution into correlated and structural (uncorrelated) contributions.

**Acknowledgments**

The authors would like to thank Dr Stefano Tarantola and Dr. Cleo Kontoravdi for their contributions and support of this work. One of the authors (O.V.K.) also gratefully acknowledges the financial support by the EPSRC grant EP/K038648/1.

**Appendix A. Sampling from joint and conditional distributions in a non-rectangular domain using cumulative distributions**

To sample from joint distribution $p(x_1,...,x_n)$ we use the following Theorem from [20]:

**Theorem 1.** Let $\gamma_1,...,\gamma_n$ be independent random numbers. The set of random values $\{\xi_1,...,\xi_n\}$ defined on $\Omega^n$ obtained from

$$\begin{aligned}F_1(\xi_1) &= \gamma_1, \\ F_2(\xi_2 \mid \xi_1) &= \gamma_2, \\ &........ \\ F_n(\xi_n \mid \xi_1,...,\xi_{n-1}) &= \gamma_n\end{aligned} \tag{56}$$

has the pdf $p^\Omega(\xi_1,...,\xi_n)$. Here $F_1(\xi_1)$, $F_2(\xi_2 \mid \xi_1)$ … $F_n(\xi_n \mid \xi_1,...,\xi_{n-1})$ are cumulative distributions corresponding to $\{\xi_1,...,\xi_n\}$.

The same approach can be used to sample from conditional distributions.

**Appendix B. Constant distribution in a triangle $\Omega_1$**

Consider an upper triangle $\Omega_1$ shown in Fig. 2 in which random variables $x_1, x_2$ have a joint pdf

$$p^{\Omega_1}(x_1, x_2) = 2, \quad (x_1, x_2) \in \Omega_1. \tag{57}$$

From the Bayes theorem:

$$p^{\Omega_1}(x_1, x_2) = p^{\Omega_1}(x_1) p^{\Omega_1}(x_2 \mid x_1). \tag{58}$$

The marginal pdf for the first variable and the conditional pdf for the second variable can be computed explicitly:

$$p_1^{\Omega_1}(x_1) = \int_{1-x_1}^{1} p^{\Omega_1}(x_1, x_2)\, dx_2 = \int_{1-x_1}^{1} 2\, dx_2 = 2x_1, \tag{59}$$

$$p_2^{\Omega_1}(x_2 \mid x_1) = \frac{p^{\Omega_1}(x_1, x_2)}{p_1^{\Omega_1}(x_1)} = \frac{1}{x_1}. \tag{60}$$

From this we can find cumulative distributions:

$$F_1(x_1) = \int_0^{x_1} 2x_1 dx_1 = x_1^2,$$

$$F_2(x_2 \mid x_1) = \int_{1-x_1}^{x_2} \frac{1}{x_1} dx = \frac{x_2 - (1-x_1)}{x_1} = \frac{x_1 + x_2 - 1}{x_1}.$$

(61)

Following (57) we use the following procedure to sample $(x_1, x_2)$:

$$x_1^2 = \gamma_1,$$

$$\frac{x_1 + x_2 - 1}{x_1} = \gamma_2.$$

(62)

After transformation

$$x_1 = \sqrt{\gamma_1},$$

$$x_2 = x_1(\gamma_2 - 1) + 1.$$

(63)

To sample $(x_1, x_2)$ and $(x_1', x_2')$ two independent sets of $(\gamma_1, \gamma_2)$ and $(\gamma_1', \gamma_2')$ are needed. To sample $(x_1', x_2 \mid x_1')$ which is needed to use formula (8) the following procedure is used:

$$x_1' = \sqrt{\gamma_1'},$$

$$x_2 = x_1'(\gamma_2 - 1) + 1.$$

(64)

This explicit sampling is limited to the case of simple geometries of $\Omega^n \subset H^n$ allowing computation of cumulative distributions. We note, that the same test case was used in [13] to illustrate the Rosenblatt transformation.